# Global virtual factory simulation for energy efficiency

Jenil Agrawal[a] | Dennis Xenos[b] | Nilay Shah[a*]

[a] *Department of Chemical Engineering, Centre for Process Systems Engineering, Imperial College London, South Kensington Campus, SW7 2AZ, United Kingdom*
[b] *Flexciton, Linen Court, N1 6AD, United Kingdom*

**ABSTRACT**

Over the past decades, industries have had to tackle the issue of sustainability as a matter of increasing urgency to mitigate greenhouse gas emissions and abide by government regulation policies. The semiconductor industry has come under scrutiny for its significant carbon and water footprints, largely attributed to substantial energy consumption. Notably, machine utilisation accounts for a major portion of this consumption, while the rest is shared by factory facilities. Thus, it is critical to optimise machine utilisation to attain desired throughput while minimising energy consumption. A flexible job shop (FJS) model is developed for the six-step Intel Minifab using a discretised time approach to find an optimal schedule. A multi-objective formulation is solved using the $\varepsilon-$constraint method to minimise energy consumption and make-span for a constant throughput. This approach affords stakeholders the flexibility in choosing machine combinations based on energy implications, a crucial consideration, especially in situations of machine unavailability. Modifying this formulation, the throughput is varied, and an additional objective is incorporated to maximise it. The energy-throughput trade-off is critically analysed, and the effect of dynamically switching off machines is well demonstrated. Energy consumed per unit wafer is observed to decrease non-linearly with make-span. The effect of optimal scheduling on energy consumption is highlighted using the FIFO heuristic strategy. When employing optimal scheduling, there was a marked reduction of 26.72% in energy consumption compared to the maximum energy consumed in the FIFO approach.

## 1 | INTRODUCTION

The fourth industrial revolution, also commonly referred to as Industry 4.0 (I4.0), gave rise to disruptive technologies and paved the path of automation in the manufacturing sector [1]. Brettel et al. [2] claim that Industry 4.0 encompasses a diverse mix of principles and tools, ranging from extensive data handling and advanced data analytics to cloud solutions, cybersecurity measures, and decentralised smart systems. These new technologies enabled higher production efficiencies covering the highly dynamic corporate environment. Traditional production systems are quite inefficient in dealing with ecological imbalances. General environmental degradation, global warming, and higher carbon footprints can all be traced back to traditional manufacturing systems [3]. Morrar et al. [4] assert that the I4.0 has the potential to mitigate the ecological and social impacts and help manufacturers contribute towards a more sustainable future. However, this practice may give rise to long-term organisational competitiveness in monetary terms.

One of the significant industries undergoing a profound transformation in the era of Industry 4.0 is semiconductor wafer production. Underpinned by substantial capital and operational expenditures, companies are constantly navigating the balance between fulfilling market demands and maintaining a competitive edge. For context, setting up a semiconductor wafer fabrication facility (fab) typically costs around 4 billion USD [5]. These fabrication facilities (fabs) are infamous for using toxic chemicals and their high energy consumption. For context, a typical medium-scale fab consumes enough electricity to power around 7000-8000 homes in the US [6]. Additionally, the semiconductor industry accounts for 1.3%-2% of the total electricity consumed by the US manufacturing sector [7]. Thus, this industry is often under pressure to achieve ecological efficiency. Generally, after the raw wafer is produced, the process is segregated into two parts: front-end and back-end [8]. The front-end production, which typically includes diffusion, etching, lithography, and implantation, requires special operating conditions in the cleanroom, leading to a large fraction of total energy consumption [9]. The cleanrooms usually consume 50 times more energy than typical commercial facilities because they require about a hundred air changes per hour done using ultra-low particulate air technique [10]. In addition to the energy demands of infrastructure, individual production machines also account for significant energy consumption. This elevates the urgency for manufacturers to give serious attention to their environmental footprint. To effectively manage this, data from individual machines is essential. However, many manufacturing units don't have adequate sensors, leading to gaps in clear information sharing [11]. This limitation muddles the finer details, making it challenging to calculate precise energy use. While companies are open about data related to their infrastructure, they often withhold specifics from the design and production stages to protect their proprietary methods.

Although stakeholders are making strides to adopt new technologies and enhance their production planning, a salient challenge of Industry 4.0 remains data collection amidst intricate and varied manufacturing processes [2]. Data on equipment and work in progress (WIP) is pivotal for stakeholders, facilitating informed production planning decisions and evaluations of environmental impacts. The academic community has been diligently refining production planning and addressing scheduling dilemmas, especially in complex sectors like semiconductor wafer fabrication, since the 1980s. This industry, in particular, stands out as a prime example of tackling these challenges, paving a seamless pathway towards the realisation of Industry 4.0 [2].

Even after these years, uncertainty in the lot (wafer) flow across the shop floor in a fab is still bemusing the scientific community. The issue of re-entrant flows and sequence-dependent setup times coupled with meeting the desired throughput is a major challenge in scheduling jobs on the shop floor. The scheduling problem can be dealt with maturity if companies prioritise their KPIs and undertake a dynamic scheduling approach. Scheduling problems have been in focus since 1956 as they directly impact production efficiency [12]. Several domains, such as manufacturing, logistics service, finance, and transportation, have always been working on incorporating optimised scheduling plans in their operations to increase their outputs. A brief overview of scheduling criteria has been summarised in Figure 1. The first criterion is the production environment. The production environment can be classified based on the types of resources available and the way they can be utilised at different production stages. In the case of single machine availability, it processes all the tasks, while in the case of parallel machines, they can carry out multiple tasks on similar multiple machines together. If a process involves multiple tasks (operations) to be carried out on multiple machines, it creates workshop conditions. The workshop can be divided into three major types based on their sequence (route) constraints: job shop, flow shop and open shop. The job shop basically encompasses the fixed sequence of tasks carried out on similar machines. While the sequence of tasks in the flow shop is not fixed and can differ based on process constraints. The open shop provides a more flexible shop floor where the sequence of tasks is not fixed, and they can be carried out in any order based on resource availability. Most of the fabs mimic job shop environments; they are characterised as flexible job shop (FJS), as the sequence of tasks might be fixed, but multiple machines have the flexibility to carry out the same task. The second criterion is based on either a dynamic or static production strategy. Dynamic strategy accounts for machine breakdowns and predictive maintenance, basically accounting for temporal changes in production. The third criterion is objective-based scheduling. This can be time-based (minimising make-span or tardiness) or job-based (maximising throughput). Other objectives used are logistics and inventory-based, which focus on downstream operations. The time-based objectives are most widely used in the academic literature to achieve company objectives. Developing a multi-objective approach considering energy consumption and factory KPIs has been catching light in recent years. Manufacturers try to develop a meta-heuristic approach to optimise energy consumption while reaching the required throughput.

The wafer fabrication shop floor presents one of the most complex scheduling problems. To understand, it has about 150 machines and a lot size of around 800 waiting to be processed. The lots flow through around 600 unique parts with different sequence-dependent setup times and multi-objective formulations, resulting in around 50,000 different setup combinations. The single test machine scheduling is a strong non-deterministic polynomial time problem ($\mathcal{NP}-$ Hard) [13]. As soon as a non-preemptive schedule is required, the problem becomes $\mathcal{NP}-$ Hard [14]. Conventional search and optimisation algorithms cannot solve these problems in polynomial time and require extensive computing resources in order to reach a solution. A common approach to solve them is using heuristic-based scheduling rules such as achieving the shortest make-span and earliest due date (EDD).

Scheduling dilemmas are typically conveyed using the conventional triplet notation $(\alpha, \beta, \gamma)$. The FJS can be symbolised as $FJ_c / M_i / C_{max}$. The initial symbol designates the shop floor type, specifically the flexible job shop. The subsequent symbol points out the machine feasibility criterion, while the final symbol denotes the chosen objective function, in this case, the make-span.

Having said that, the structure of this paper is organised as follows: Section 2 provides a brief description of academic literature encompassing energy-efficient simulation and energy calculations. It is followed by describing the Intel Minifab model, the concept of state task networks and formulating a time discretised FJSS problem (FJSSP) in section 3. The benchmark test is performed and analysed in section 4. Section 5 focuses on incorporating energy consumption parameters, highlighting the essence of scheduling. The paper concludes with an outlook in section 6.

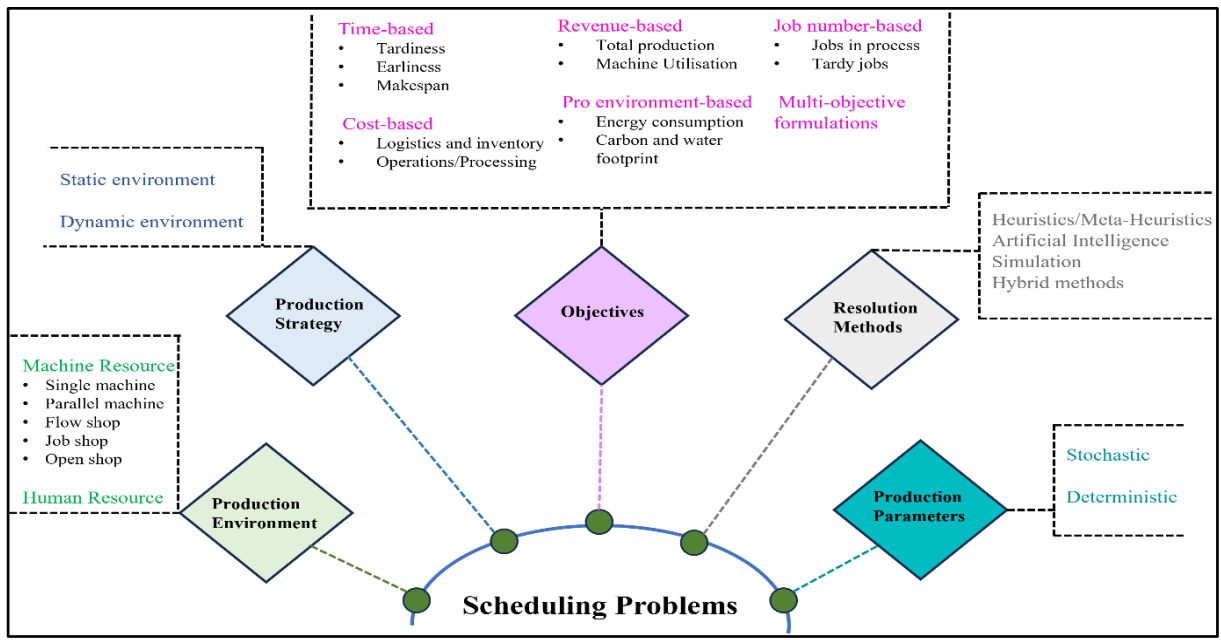

**Figure 1.** Scheduling problems characterisation

## 2 | RESEARCH BACKGROUND

Extensive research, encompassing methodologies from simulation-based approaches to parametric analyses, has been conducted to analyse energy consumption within the fabs and to elucidate the associated carbon and water footprints. An energy-aware model based on queuing theory considering re-entrant flows was developed by Hyun et al. [15]. The model focuses on operating idle machines in low energy state to investigate optimal energy consumption. Rahimifard et al. [16] develops a plant-level model to break down energy consumption in a top-down approach starting from plant to process level. They analyse energy as fixed (facilities) and variable (machine utilisation) consumption. An alternate approach using an operating curve and queuing theory, a discrete event simulation (DES) model, is developed by Schneider et al. [17]. Main focus is laid on resource utilisation v/s flow factor (ratio of cycle time and raw process time). Later, the flow factor is also compared against the ecological factor. This study is by far the most relevant to comment on how resource utilisation affects energy consumption. Nouiri et al. [18] conducted an extensive review of modern energy-conscious dynamic scheduling methods used in production and logistics systems, including inventory management. They proposed a flexible, dynamic scheduling model that considers key performance indicators, like how well and efficiently a facility operates. Their findings were confirmed through experiments. However, their models mainly focus on machine breakdowns and don't consider other aspects like re-entrant flows that might require job rescheduling. Predominantly, academia has approached scheduling models without explicit consideration of energy dynamics, opting instead to derive energy implications from the resultant resource allocation. A subset of this literature also discusses the energy expenditure of the manufacturing infrastructure and non-production-centric resources. A critical review by Mullen et al. [19] underscores the viability of a cost-efficient alternative lithography method, specifically block-copolymer, as a solution to address the challenges associated with miniaturisation. A simulation study by Kircher et al. [20] demonstrates that using heat recovery systems for exhaust air yields an 11.4% reduction in total energy consumption. Wang et al. [21] conducted an empirical analysis wherein they ascertained that the clear dry air (CDA) system was responsible for 19.8% of the energy utilisation among the evaluated systems. Furthermore, the researchers posited that through the implementation of exhaust heat recuperation techniques, coupled with the mitigation of heat losses and rectification of pipe leakages, there is a prospective reduction in energy consumption by approximately 2.4%. Despite the extensive research on energy consumption within manufacturing facilities, a critical question still persists: How does the selection of machines through an optimal schedule influence the overall energy footprint of a fab?

In contrast to many studies being reported, this paper focuses on the front-end process in a fab and uses a flexible job shop scheduling (FJSS) formulation of a six-step-Minifab model proposed by Dr. Kempf (associated with Intel) [22] to incorporate energy analyses based on machine utilisation while attaining the maximum throughput in the chosen time horizon. To the authors' knowledge, no such study has been published focusing on complex semiconductor wafer fabrication facility.

## 3 | MODEL DESCRIPTION

### 3.1 | Intel Minifab Model

Developing a full-scale semiconductor manufacturing model is tedious and erroneous. Thus, researchers use available testbeds and datasets to simulate their models. There are four widely used datasets in academia: MIMAC datasets, Minifab, Harris (Kayton), and SEMATECH 300 mm. Their properties are briefly described in Table 1.

**Table 1.** Testbed/Dataset properties

| Dataset | #Machines/ Workstation | #Products | #Process steps |
|---|---|---|---|
| **MIMAC** [23] | <260 / <85 | <21 | <280 |
| **Harris** [24] | 12 / 11 | 3 | <22 |
| **Intel Minifab** [25] | 5 / 3 | 2 | 6 |
| **SEMATECH** [26] | 275 / 103 | 1 | 364 |

The Minifab model proposed by Kempf and Spier [25] is a scaled-down model that accurately captures actual fab characteristics. The model contains re-entrant flows, sequence-dependent setup times, batch process and operator assignments. It has three workstations, featuring six steps with five machines. The process flow in the Minifab model is shown in Figure 2. The three workstations are diffusion, implantation, and lithography. Since there are three workstations and each wafer requires six steps (operations), each wafer goes through each station twice before coming out as a finished product. Process lines 3,4, and 5 are the re-entrant lines which make it possible. Each diffusion and implantation have two eligible machines each – machines A and B and machines C and D, respectively, while implantation has only one machine, machine E. The implantation centre is the bottleneck workstation as it only has one eligible machine. Machines A and B are identical and require 75 mins of predictive maintenance per day, while machines C and D are identical, requiring 120 mins of predictive maintenance per 12-hour shift. There are two types of products – A and B, starting with a volume of 50 and 30 lots per week. The model is based on a cellular layout with an automated material handling system (AHMS). The AHMS can only handle one lot at a time, and hence, operator assignments are introduced to transfer multiple lots across one or more workstations.

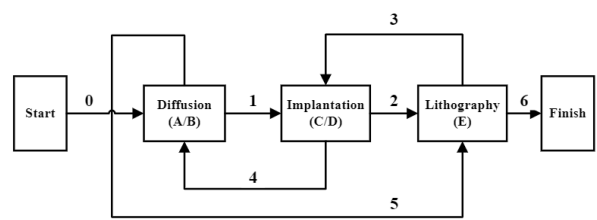

**Figure 2.** Intel Minifab- process flow representing all workstations with process lines. The process lines are numbered to signify the order of operation.

The objective is to determine an optimal schedule accommodating the bottleneck machine, machine E, which is giving sub-optimal fab performance. This study focuses on developing a FJSS formulation and lays an additional focus on analysing energy consumption based on machine utilisation. Hence, the Minifab model is simplified further, and the following assumptions have been made:

1. No travel time is modelled. It is assumed the transferring of wafers from one station to another doesn't require any additional time.
2. No operators or scheduled maintenance is considered. Neither do we consider any loading or unloading time.
3. Processing time for each operation is deterministic.
4. All operations are non-preemptive (an operation, once started, cannot be interrupted/stopped before completion)
5. No capacity limits considered.
6. No batching of lots or buffers is considered. Each individual wafer is treated as a single batch.
7. Process time unit is hours.

### 3.2 | State Task Networks

As the complex flow patterns (re-entrant lines in this case) give rise to certain ambiguities in traditional flowsheet formulation, we use the concept of the state-task network (STN). STN was introduced by Kondili et al. [27] to formulate a short-term scheduling algorithm for batch processes. The idea is to break the complex process lines into two types of nodes - states and tasks, to simplify the network. The state nodes represent the feeds, intermediate and final products, while the tasks represent process (operation) steps which transform one state to another or gives rise to multiple states based on the type of operation and input state. States and tasks are denoted by circles and rectangles, respectively, connected by directional arrows. The basic rules for constructing a STN are:

1. If two or more streams enter the same task, they should be of the same quality (no mixing allowed in the states). If the process involves mixing, a separate task should be created for that operation.
2. The number of input (output) states in a task depends on the number of types of material/product present in those input (output) states.

It should be noted that STNs can be disjoint graphs. For example, there are certain manufacturing processes where products don't share the same resources (same all machines) but are coupled as they use the same manufacturing facilities and might be sharing some common resources. The STN for such processes are not necessarily connected and can be formulated as disjoint sub-graphs. A detailed example describing the STN for a multipurpose batch plant, which considers sequencing as well as scheduling simultaneously, can be found in [27].

For the purpose of this paper, we are interested in building STN for the Minifab model described in section 3.1. As we know, each wafer goes through each workstation twice, it allows us to separate each workstation into two identical stations (tasks). For example, diffusion to diffusion 1 and diffusion 2; similarly, we get implantation 1, implantation 2, lithography 1 and lithography 2. This implies that now each wafer will go through these six stations (tasks), which might be identical pairs but produce different states. For example, after completing the operation at each workstation once, the wafer re-entering the diffusion station again would not be the same as the new wafers entering the same diffusion station, and thus, they cannot be allowed to enter process line 1. Hence, we segregate each workstation into two individual stations. It is interesting to note that the re-entrant lines are characterised by machines' availability. In simpler terms, if a wafer returns to the same station after being processed once and the station has limited machines to process new wafers as well as re-entered wafers, the assignment of the wafer to an available machine is what gives rise to complexity in scheduling re-entrant flows. Hence, utilising this concept, our STN captures the re-entrant flow as each pair of tasks will havethe same eligible machines as described in the original model. For example, task 1 - diffusion 1 can still be processed by machines A and B and simultaneously, the same machines will be shared with task 2 - diffusion 2. The state task network for the Minifab model is represented in Figure 3.

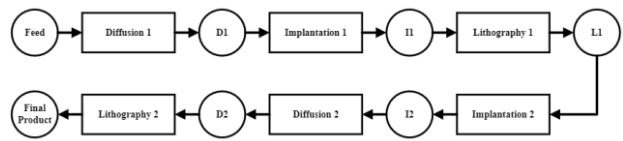

**Figure 3.** STN for the Minifab model.

### 3.3 | Mathematical formulation for the scheduling problem

In this paper, we formulate a mixed integer linear programming (MILP) model to solve the scheduling problem. The mathematical model ensures the listing of key constraint sets to capture operation precedence and feasible operation assignment to eligible machines by choosing an appropriate set of parameters and variables.

One of the key considerations in any mathematical formulation is the choice of time representation. In our study we use a discretised time approach. The time horizon is divided into a number of equal time intervals. It is evident that start or end of any operation should occur within these interval boundaries. We discretise the time intervals from 1 to T, starting at 1, and final time would be T+1. Each interval is taken as 1 time unit. This approach provides a temporal grid to track the operations and machine utilisation and it has been used by scholars previously to demonstrate production planning and scheduling [28].

| Parameters | |
|---|---|
| $j$ | Represents jobs |
| $k$ | Represents operations |
| $m$ | Represents machines |
| $t$ | Represents time periods |
| $p_{k,m}$ | Represents process time for $k^{th}$ operation on $m^{th}$ machine |
| $S_{k,k'}$ | Represents operation sequence |
| | |
| Variables | |
| $X_{k,j,m,t}$ | Binary variable taking value 1 if operation $k$ of job $j$ starts on machine $m$ at time $t$ |
| $W_{k,m,t}$ | Binary variable taking value 1 if machine $m$ is available to carry out operation $k$ at time $t$ |
| $C_{k,j,m,t}$ | Continuous variable to capture completion time $t$ for operation $k$ of job $j$ on machine $m$ |

*Constraints*

1. Operation initiation constraint

$$\sum_t \sum_m X_{k,j,m,t} = 1 \qquad \forall j, k \qquad (1)$$

The operation initiation constraint ensures that every operation associated with each job is initiated precisely once. This is crucial for maintaining the integrity of the production schedule, preventing operation duplications, and ensuring that each job undergoes every required processing step.

2. Machine availability/assignment constraint

$$\sum_j \sum_{t'=t}^{t-p_{k,m}+1} X_{k,j,m,t} \leq 1 \qquad \forall k, m, t \qquad (2)$$

This constraint ensures, operation $k$ can be performed at the machine $m$ at time $t$ if no other operation is performed on the same machine during its execution. It is interesting to note that we implement a 'full backward aggregation' approach for machine allocation constraint as described by Shah et al. [29]. Constraint (2) 'looks backward' in time and check for potential clashes in resource allocation. Although it might not impose tightest of bounds, it is advantageous in reducing the size of our MILP formulation, as the number of constraints used in this case are lower as compared to the 'forward aggregation' approach which looks forward in time $(t + p_{k,m} - 1)$. More details regarding the forward aggregation approach can be found in [29]

3. Machine-operation eligibility constraint

$$X_{k,j,m,t} \leq 1\{p_{k,m} > 0\} \qquad \forall k, j, m, t \qquad (3)$$

In the above constraint, $1\{\cdot\}$ is an indicator operator, which is 1 if the condition is true and 0 otherwise. Constraint (3) makes sure that the job $j$ doesn't start its operation $k$ unless it is assigned an eligible machine $m$ for that operation. This is done by checking the non-zero values in $p_{k,m}$ parameter because if the processing time for any operation on a machine is zero it is deemed to be in-eligible.

4. Parallel processing constraint

$$X_{k,j,m,t} + X_{k,j',m',t} \leq 2 \qquad \forall k, j, m, t \qquad (4)$$

This constraint allows operating two jobs simultaneously on two eligible machines if they are available.

5. Sequential processing constraint

$$\sum_{m'} \sum_{t'=t}^{t+p_{k,m}-1} X_{k',j,m',t'} \leq X_{k,j,m,t} \qquad \forall k, j, m, t, S_{k,k'} \qquad (5)$$

Constraint (5) suggests that the job cannot be processed simultaneously on two machines. It also implies that succeeding operation $k'$ can start only when preceding operation $k$ is executed on any eligible machine on the shop floor. It should be noted that each job follows a fixed sequence of operations, and hence we don't use big-M method to decide precedence among the operations.

6. Capturing completion times of operations

$$C_{k,j,m,t} = X_{k,j,m,t} * (t + p_{k,m} - 1) \qquad (6)$$

The completion times of operation $k$ for job $j$ ending at time period $t$ is calculated by constraint (6).

7. Releasing new job constraint

$$X_{k,j,m,t} \leq X_{k,j-1,m,t-p_{k,m}} + \sum_{m'} W_{k,m',t}$$
$$\forall\ k,j,m,t \quad (7)$$

This constraint mandates that a job $j$ can commence its operation $k$ on machine $m$ at time $t$ only if its predecessor, job $j-1$, initiated the same operation $k$ on the same machine $m$ at an earlier timeframe or if alternate machines are accessible to execute the job at the given time. This formulation ensures an orderly and efficient processing of products, respecting their sequence and the machine availability.

While at first glance, constraint (5) and (7) appear to address similar scheduling aspects, their roles within the model are distinct. Constraint (7) governs the release of a job based on the start time of its predecessor, thereby offering enhanced flexibility in the job flow. In contrast, constraint (5) ensures that a job is only released once the preceding job has completed its designated operation. Though each of these constraints independently may not impose stringent scheduling conditions, their combined application ensures a robust scheduling framework by imposing tighter bounds.

8. Objective function

In the range of time-based objective functions discussed in section 1, we choose to minimise 'make-span' ($C_{max}$) for our model. The strength of this objective lies in its ability to both minimise job waiting times and ensure a balanced distribution of work across resources, preventing any particular machine from being consistently exploited.

$$Min\ C_{max} \quad (8)$$

Additionally, constraint (9) captures the make-span by considering the completion times for the last operation for all jobs.

$$C_{max} \geq C_{k,j,m,t} \quad \forall\ j \quad (9)$$

9. Constraints demonstrating nature of decision variables and choice of time horizon.

$$X_{k,j,m,t}, W_{k,m,t} \in \{0,1\} \quad (10)$$
$$C_{k,j,m,t} \geq 0 \quad (11)$$

## 4 | BENCHMARK TEST

This section discusses a benchmark test for the MILP formulation discussed in section 3.3. The process times for each operation at its eligible machine is shown in Table 2. It should be noted that term 'wafer' is used instead of 'job' for this section to capture the essence of discussion.

**Table 2.** Process times for operations in the Minifab model

| Operation | Eligible Machines | Process time(hour) |
|---|---|---|
| Diffusion 1 | Diffuser 1, Diffuser 2 | 2 |
| Implantation 1 | Implanter 1, Implanter 2 | 1 |
| Lithography 1 | Lithographer | 2 |
| Implantation 2 | Implanter 1, Implanter 2 | 1 |
| Diffusion 2 | Diffuser 1, Diffuser 2 | 1 |
| Lithography 2 | Lithographer | 2 |

It should be noted that the values for process time chosen for this study are scaled down values in hours without the loss of generality. The original values can be found in [25]. The process time in the original model is described in minutes and given our discretised time approach, larger process time values that results in choosing longer time horizon. This makes the model computationally expensive and inefficient.

The code was developed using General Algebraic Modeling System (GAMS)[30] and solved for 5 jobs using CPLEX solver. MILP is solved in CPLEX using a branch and bound algorithm [31]. The resultant minimum make-span is **17 hours**. The schedule for each job (wafer) is presented in the form of Gantt charts in Figure 4. All the wafers are considered similar, and no precedence priority is assigned to a wafer to begin its operation. We can observe wafer 1 and 3 start diffusion 1 parallelly on available diffuser 1 and diffuser 2 at t = 1, followed by wafer 2 and wafer 4. It is noteworthy that wafers 1 and 2 are processed without any waiting times, while other wafers experience a delay of at least 3 hours before starting lithography 1. This can be attributed to a single eligible machine available for lithography, and thus, it can also be regarded as a bottleneck workstation. A significant observation is a one-hour delay before initiating implantation 1 (from t = 3 to t = 4) for wafer 3, even when only wafer 1 is undergoing implantation 1 at that interval. This suggests that an implanter should be available during this period. However, the optimisation algorithm opts for this delay because, even if it processes wafer 3 during that window, the wafer will inevitably encounter a wait time before lithography 1. Thus, exploiting additional machines under such circumstances wouldn't yield any tangible result in terms of make-span. Additionally, as the wafers progress through intermediate stages, the inherent variability in processing times starts playing a significant role. Some wafers might spend more time in one process compared to others. Consequently, the originally clustered group of wafers begin to disperse, leading to a more spaced-out flow, and hence we don't observe any waiting times before implantation 2 and lithography 2.

Increasing the quantity of wafers invariably leads to extended waiting periods within the production facility due to limited machinery resources. To check this, we simulate the model for 15 wafers. The minimum make-span obtained is **37 hours**. The Gantt charts for 5th wafer and 15th wafer is shown in Figure 5. We observe that the total operational time for wafer 5 remains consistent at 17 hours, akin to prior observations. However, a consistent delay of 2 hours emerges prior to both implantation 1 and lithography 1, in contrast to the previous 1-hour and 3-hour intervals. This shift arises as the solver anticipates more wafers, which are yet to be released from the feed and aims to mitigate congestion before lithography 1. As predicted, the waiting time prior to implantation 1 for wafer 15 is increased to 3 hours.

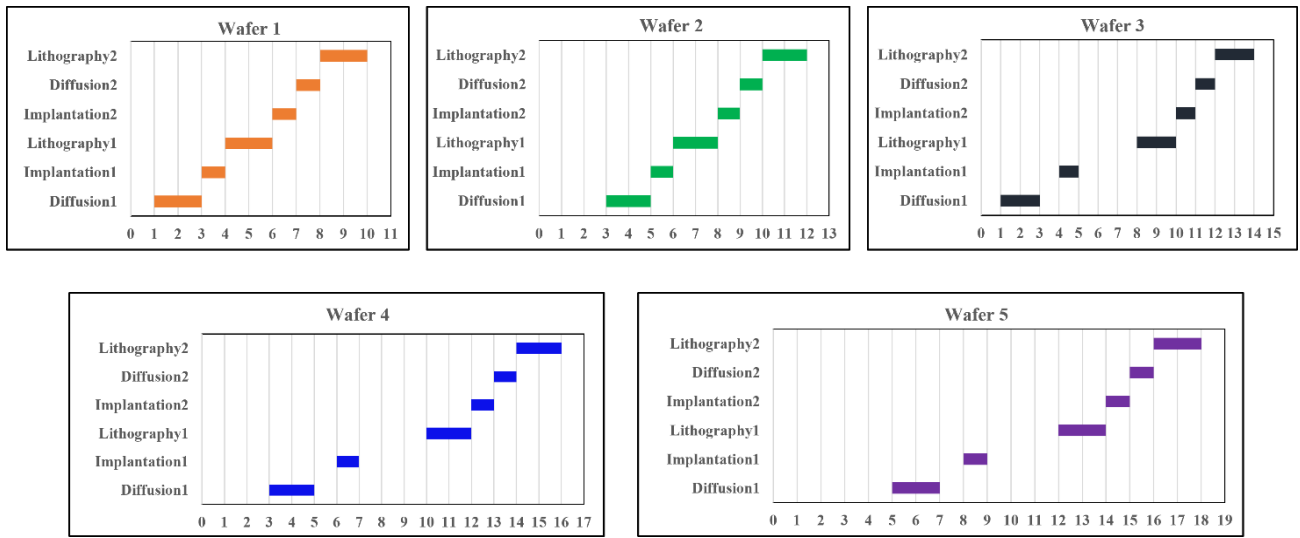

**Figure 4.** Gantt charts representing optimal schedule for each wafer in the test case with total 5 wafers.

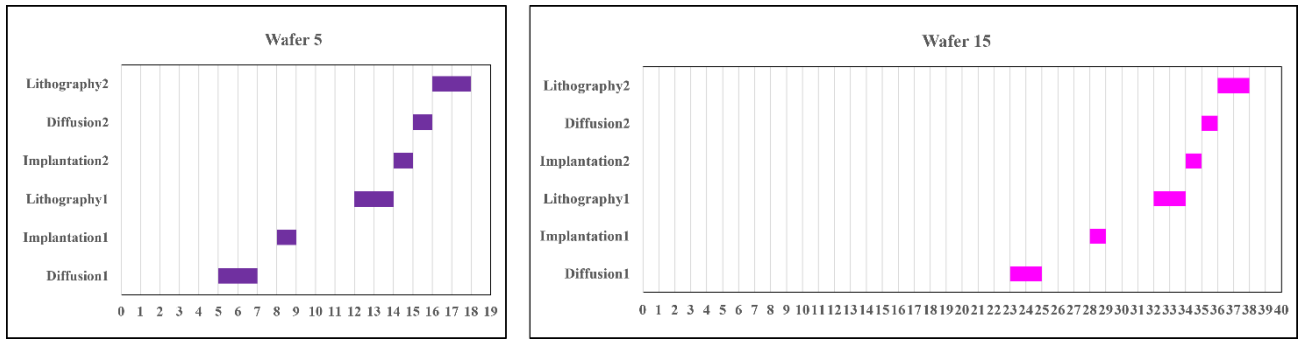

**Figure 5.** Gantt charts representing optimal schedule for wafer 5 (purple) and wafer 15 (pink) in the test case with total 15 wafers.

## 5 | ENERGY CONSUMPTION MODEL

This section focuses on understanding and incorporating energy consumption parameters in the FJSS model described in section 3. As discussed in section 1, fab is an energy-intensive shop floor. We can break the energy consumption into two types: fixed and variable. A pictorial representation of breakdown is shown in Figure 6. The fixed energy consumption can be attributed to the factory operating conditions such as clean room maintenance, operating heating ventilation and air conditioning (HVAC) systems, and chillers which makes the major fraction of the total energy consumption [20]. This consumption is regardless of volume of production and depends on the time for which fab is operating or can also be understood as a function of make-span. The operating machines (switched on) even though producing/not producing can be accounted for variable energy consumption. This is variable because it depends on number of machines being used in the given time, if a machine is switched off it won't be consuming energy. Hence, variable energy consumption also affects the throughput. A natural instinct would be to exploit all the machines and achieve the maximum throughput. But in doing so we are also increasing the share for fixed energy consumption which is much higher. Thus, it is critical to analyse the trade-off between energy consumption and throughput for a given time horizon.

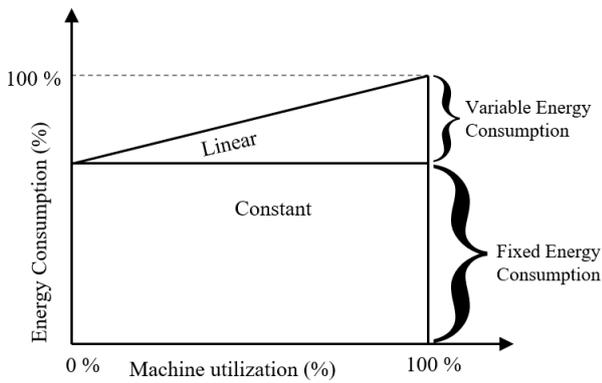

**Figure 6.** Concept curve to demonstrate fixed and variable energy consumption.

## 5.1 | Energy consumption analysis with fixed throughput

In this section we lay a focus on variable energy consumption based on machine choice and utilisation. Machines are assigned weights signifying power rating and two additional equations are introduced for energy consumption calculations.

Additional equations :

$$E_m = \sum_t \sum_k \sum_j X_{k,j,m,t} * p_{k,m} * wt_m \quad \forall\, m \quad (12)$$

$$TE = \sum_m E_m \quad (13)$$

here $wt_m$ represents power rating (weights) of each machine. The equation (12) calculates energy consumption for each machine based on its utilisation, while equation (13) accounts for total energy consumption and acts as a second objective function in our formulation. Now we have a multi-objective formulation to minimise make-span as well as energy consumption for a given throughput and is solved using the $\varepsilon-$ constraint method [32]. At an initial glance, the weighted-sum method might appear appealing and direct, as it forms a single objective using a linear combination, especially when we aim to minimize both objectives. However, the $\varepsilon-$ constraint method is favored for generating the Pareto front, given its robustness in handling non-convex Pareto fronts, a scenario where the weighted-sum approach might fall short.

The model is run for 15 wafers and the power rating for machines are shown in Table 3. It should be noted that these are dummy values chosen to understand energy consumption across the shop floor.

**Table 3.** Power rating (weights) of machines

| Machine | Power rating (KW) |
|---|---|
| Diffuser 1 | 0.01 |
| Diffuser 2 | 0.02 |
| Implanter 1 | 0.1 |
| Implanter 2 | 0.2 |
| Lithographer | 0.1 |

The bounds of Pareto front are obtained by minimising each objective function. The energy consumption is discretised linearly withing these bounds and each discretised value act as an upper bound (constraint) for $TE$ and we minimise make-span as primary objective.

Interestingly, the approach gives us distinct energy consumption based on different machine utilisation, for the same minimum make-span. Hence for each operating condition (choice of machines) we get varied energy consumption. The algorithm chooses optimal set of machines to achieve the throughput (15 wafer) with minimum make-span, thus giving us a flexibility in choosing machine combinations based on energy usage. The minimum total energy obtained is 9.450 KWh while maximum is 10.190 KWh. Set of machine choices corresponding to minimum and maximum energy consumption is reported in Table 4. It is evident that we can achieve the desired throughput utilising only three machines with minimum total energy. But it is not always the case in a real factory as there might be instances when one or many of these machines are unavailable due to maintenance or breakdown and hence flexibility in choosing resources is quite necessary.

**Table 4.** Machine utilisation instances for min v/s max energy consumption

| Machine | Time utilized (hour) (Min Energy) | Time utilized (hour) (Max Energy) |
|---|---|---|
| Diffuser 1 | 45 | 31 |
| Diffuser 2 | 0 | 14 |
| Implanter 1 | 30 | 24 |
| Implanter 2 | 0 | 6 |
| Lithographer | 60 | 60 |

## 5.2 | Energy consumption analysis with variable throughput

In a fab, a high throughput often guarantees higher revenue, but it comes at the cost of energy usage. With the energy prices rising it would be critical to analyse the trade-off between throughput v/s energy. One idea would be considering energy aware scheduling based on 'time of use' (TOU) pricing. A study by Park et al. [33] focuses on FJSSP based on TOU and scheduled maintenance. Their results demonstrate a 6.9% reduction in energy cost while maintaining maximum productivity. But they do not control machine processing times with the TOU pricing simultaneously. Other idea would be varying the make-span, and for given make-span we maximise the throughput and minimise the energy consumption. This will answer an important question- how does the energy consumption vary marginally with increasing throughput? This section employs second approach, and the following changes are made in the existing model,

1. Equation (1) is modified to,

$$\sum_t \sum_m X_{k,j,m,t} = y_j, \quad \forall\, j, k \quad (14)$$

where $y_j$ is a binary variable taking value 1 if a job j is being processed and 0 otherwise.

2. An additional equation to calculate throughput is introduced,

$$Tp = \sum_j y_j \quad (15)$$

Additionally, now the model has an additional objective to maximise the throughput. The maximum make-span chosen here is 80 hours with a maximum throughput of 50 wafers. As discussed $\varepsilon -$constraint method is used, make-span and energy consumption are linearly discretised to form upper bounds for respective objectives in each iteration while the primary objective is to maximise the throughput.

Figure (7) shows the throughput - make-span - energy consumption interplay. The throughput increases with increasing make-span, which is quite trivial, consequently energy consumption increases. A maximum of 37 wafers can be fabricated in a make-span of 81 hours. It is crucial to note that initially up to 8 hours, throughput is 0, as it would take minimum of 9 hours to complete all operations on a single wafer. Energy consumption during this time is also 0 as the machines don't operate on any wafer during this time interval. While for maximum throughput it is 26.81 KWh. The gradient of throughput with make-span is 0.456, while for energy with make-span is 0.330. This means rate of throughput production increases at a higher rate than energy consumption with make-span. This result validates stakeholder's argument that running a fab for longer durations yields more revenue as revenue from a wafer exceeds energy costs at longer time horizons. It also points towards an argument that the energy consumed per unit wafer decreases with increasing make-span.

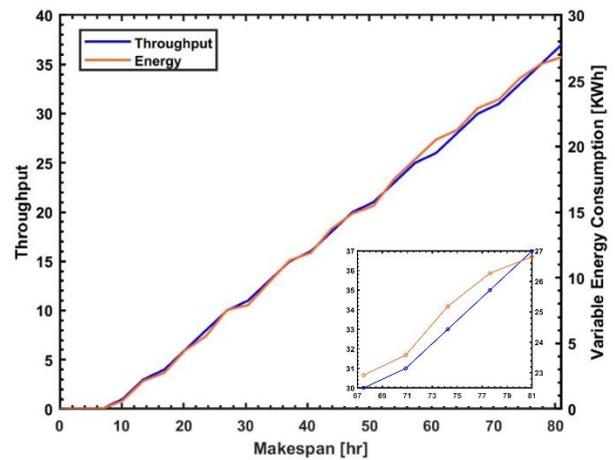

**Figure 8.** Throughput v/s Make-span plot with corresponding variable energy consumption.

Figure 8 plots energy consumption against throughput to understand the trade-off. While the plot shows a linear correlation between the two, a close observation yields that gradient of energy with throughput decreases after a throughput of 31. Thus, it again signifies, fabricating more wafers will bargain the energy usage after certain threshold value.

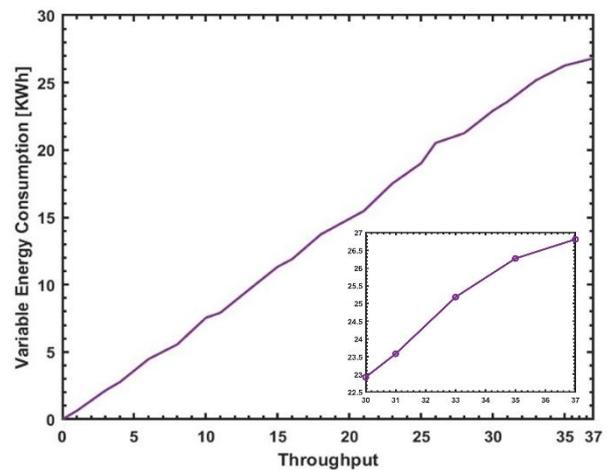

**Figure 7.** Energy v/s Throughput trade off curve.

A salient question stemming from the aforementioned results is: given an extended time horizon to achieve a consistent throughput, say 37, would the algorithm leverage this added temporal flexibility to judiciously select machines, thereby optimizing energy consumption? This can be answered using a short sensitivity analysis. The following objective is used for a fixed throughput of 37 wafers,

$$Objective = \alpha * makespan + (1 - \alpha) * TE \quad (16)$$

where $\alpha$ is a weight parameter ranging from 0 to 1.

The time horizon is extended to 120 hours, while throughput is fixed at 37. When $\alpha =1$, prioritizing minimum make-span, the solver delivered a make-span of 81 units with an energy consumption of 26.81 units, consistent with previous results. While for intermediate values of $\alpha$, the solver gives an energy consumption of 23.31 KWh with a make-span of 81 hours. Peculiarly, as $\alpha$ transitions from 0.9999999 to 1, energy consumption very gradually rises, culminating at 26.10 KWh for $\alpha =1$ as solver assigns more weight to minimise make-span with increasing alpha. On the other hand, even for a very small $\alpha$ value of 0.00000001, emphasizing energy minimization, the make-span remained unchanged at 81 units with an energy of 23.31 KWh. Even increasing the $\alpha$ very gradually results in minimum make-span of 81 hours. However, decreasing α further results in the solver prioritizing energy reduction, leading to a utilization of the entire time horizon and a resultant make-span of 120 (max time horizon). In essence, while energy displays some sort of sensitivity to α values closer to 1, make-span remains largely unaffected across the range, except when α (~ 0 ) is reduced substantially.

### 5.3 | Dynamic machine utilisation-based energy consumption analysis

So far energy consumption was calculated based on machine process times, assuming they remain on only for the time being utilised. But practically this is not the case, as a machine might be switched on, remain in an idle state, and not process any wafer. Moreover, machines should not be switched on/off frequently as this can lead to machine breakdown and hence an adaptive model must be built to account for wisely switching on/off machines. This section focuses on formulating an adaptive machine utilisation model for optimal scheduling. The following assumptions are introduced:

1. All machines are switched on at t=1.
2. Each machine has an "initialization time" required after being switched on before it can begin processing tasks. Crucially, this initialization duration is incurred every time the machine restarts, rather than being a one-off cost.
3. Once machine is switched on it must remain on for a specific minimum 'on time' before it can be turned off.

Mathematically the following changes are made in the model discussed in section 3.3 and 5.2,

| Parameters | |
|---|---|
| $st_m$ | Represents initialisation time (start-up time) for machine $m$ |
| $z_m$ | Represents minimum time for which machine $m$ must be kept on |
| | |
| Variables | |
| $o_{m,t}$ | Binary variable taking value 1 if machine $m$ is on at time $t$ |
| $u_{m,t}$ | Binary variable taking value 1 if machine $m$ is switched on at time $t$ |
| $v_{m,t}$ | Binary variable taking value 1 if machine $m$ is switched off at time $t$ |

Constraints

1. Modified machine-operation eligibility constraint

As we consider a 'startup time' after which machine $m$ can begin its operation $k$ once it is switched on, we modify equation (3) to a accommodate this criterion as below,

$$X_{k,j,m,t} \leq 1\{p_{k,m} > 0\} * (1 - \sum_{t'=t-st_m}^{t} u_{m,t}) \quad (17)$$

2. Machine is on while processing a job

This constraint ensures that if a machine m is carrying operation $k$ on any job $j$ at time $t$, then the machine must be in "on" state at that time. Essentially, it guarantees that a machine cannot process a job unless it is switched on.

$$o_{m,t} \geq \sum_k \sum_j X_{k,j,m,t} \quad (18)$$

3. Machine switch on event

The constraint ensures that if the machine is transitioning from an OFF state to an ON state, there must be a corresponding switch-on event. This constraint ensures the model correctly captures the moments when a machine is turned on.

$$o_{m,t} - o_{m,t-1} \leq u_{m,t} \quad (19)$$

4. Machine switch off event

The constraint ensures that if the machine is transitioning from an ON state to an OFF state, there must be a corresponding switch-off event. In simpler terms, this constraint ensures the model correctly captures the moments when a machine is turned off.

$$o_{m,t-1} - o_{m,t} \leq v_{m,t} \quad (20)$$

5. Switch on criteria

This constraint ensures that a machine is instructed to switch on only if it was in an OFF state in the preceding time period. It safeguards against illogical scenarios where a machine, already in the ON state, receives multiple consecutive turn-on instructions.

$$u_{m,t} + o_{m,t-1} \leq 1 \quad (21)$$

6. Mutually exclusive switching constraint

This constraint ensures machine m is not switched on and off at the same time $t$. Both events are mutually exclusive.

$$u_{m,t} + v_{m,t} \leq 1 \quad (22)$$

7. Updating machine status

This constraint updates if machine is on/off based on previous switching operations. Essentially if machine is off at time $t-1$ and it was switched on at time $t$, machine status will be changed to 'on'. Similar is the case in switching off the machine.

$$o_{m,t} = o_{m,t-1} + u_{m,t} - v_{m,t} \quad (23)$$

8. Prevent switching off during operation

A machine cannot be turned off while it is carrying out an operation, but it can do so once the operation is complete if required.

$$\sum_{t'=t}^{t' = t+p_{k,m}} v_{m,t'} \leq 1 - X_{k,j,m,t} \quad (24)$$

9. Minimum on time constraint

This constraint ensures that the machine $m$ needs to remain switched on for a certain minimum time before it can be switched off if required.

$$\sum_{t'=t-z_t+1}^{t'=t} o_{m,t'} \geq z_m * v_{m,t+1} \quad (25)$$

### 5.3.1 | Results and Discussion

The model was formulated as a multi-objective problem, as discussed in section 5.2. The make-span was varied, to maximise throughput and minimise the energy consumption using an $\varepsilon$ −constraint discretising make-span and energy. Again, the maximum time horizon chosen is 80 hours with maximum throughput of 50. The 'initialisation time' $st_m$ and 'minimum on time' $z_m$ values are described in Table (5) and (6), respectively. These values are chosen as test values to understand the fab behaviour mimicking a more realistic scenario.

**Table 5.** Time taken by machine to start (start-up time)

| Machines | Startup time (hour) |
|---|---|
| Diffuser 1 | 1 |
| Diffuser 2 | 1 |
| Implanter 1 | 2 |
| Implanter 2 | 2 |
| Lithographer | 1 |

**Table 6.** Minimum run time of the machine

| Machines | Min time (hour) |
|---|---|
| Diffuser 1 | 6 |
| Diffuser 2 | 6 |
| Implanter 1 | 4 |
| Implanter 2 | 4 |
| Lithographer | 8 |

Figure 9. shows variation of throughput and energy consumption with make-span. It is observed that the throughput and energy consumption increase linearly with make-span. But it is interesting to note that even though we get our first wafer produced at the t~10$^{th}$ hour, some energy consumption is observed. This is because there might be certain intervals <t=8 when the machine might be switched on, but it was not able to complete the operation. Furthermore, the changing gradient of energy consumption with respect to throughput offers insights into the efficiency of the production process. While the initial stages see energy consumption rise at a steeper rate than throughput, a pivot point is reached around the throughput threshold of approximately 20 units. Post this threshold, the efficiency gains become more pronounced, with the rate of throughput significantly outpacing the rate of energy consumption. This inflection point indicates there exists an optimal throughput value after which energy costs are outperformed by production revenue. The diminishing rate of energy consumption, especially noticeable for throughput values greater than 30, can be attributed to more consistent machine operations. Machines that remain operational for extended periods can bypass the repeated 'startup times', thereby achieving a state of continuous production. Each startup phase, though seemingly negligible, contributes to energy expenditure without corresponding production output. By reducing the frequency of these start-ups, the system not only saves energy but also achieves a higher operational efficiency.

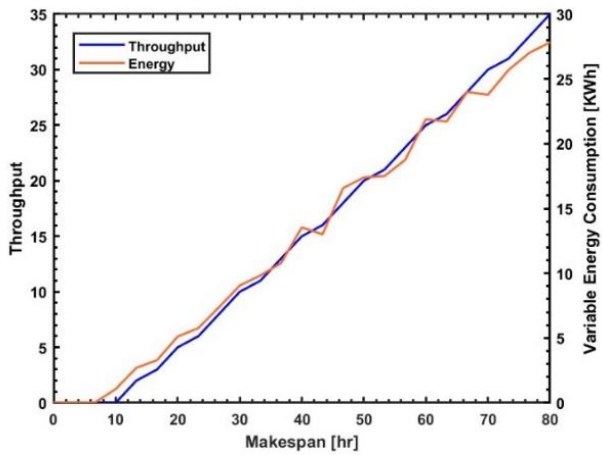

**Figure 9.** Throughput and corresponding energy consumption for dynamic machine switching model.

Another important observation is, in this case, we obtain a throughput of 35 wafers, which is two wafers less than what we obtained from the 'baseline' model described in section 5.2. This is a manifestation of 'startup time' introduced to mimic a more practical shop floor. Stakeholders would ideally like to keep the machines on and operating at all times to get maximum production output, but it is not possible due to maintenance, breakdown, and limited human resources. Although on the brighter side, if we check the energy consumption, it is 27.8 KWh. This is the energy consumption accounting for the time when machines are on, even though they are operating or non-operating. If we calculate the energy consumption in a similar fashion (time for which the machine is on) for the baseline model, it is 34.9 KWh. So essentially, two fewer wafers are fabricated, but a ~7KWh reduction in energy consumption is observed. This is the case because the algorithm switches off machines which might be idle for a longer duration and might give a sub-optimal schedule. The strategy of dynamically controlling machine states leads to energy savings. Additionally, there is a heteroscedastic increase of difference in energy consumption among both models. This trade-off can be used by small and mid-sized fabs to incorporate adaptive switching. In the long run, this approach could also extend the lifespan of the machines by preventing unnecessary wear and tear, leading to reduced maintenance costs and longer intervals between machine replacements. Moreover, even though we don't consider human resources as a parameter in this study, the abovementioned approach would allow managers to allocate human resources more effectively. Idle (switched off) machines don't need monitoring, allowing staff to focus on other critical areas of the operation.

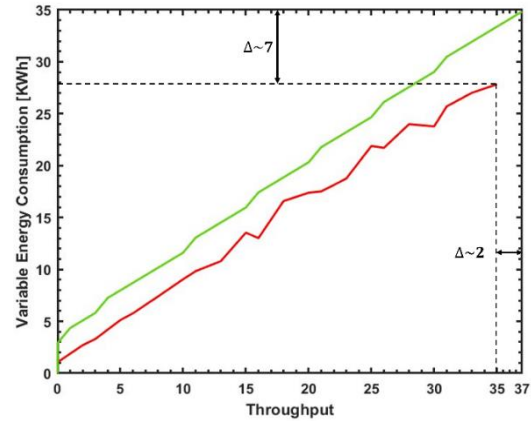

**Figure 10.** Comparison of energy consumption and throughput for baseline model (green solid line) and dynamic machine switching model (red solid line).

One more speculation, given the above observation, is that energy consumed per unit wafer will decrease. This can be validated from Figure 11. The plot is constructed from the time when at least 2 wafers are fabricated by each model. A non-linear decrease in energy consumed per unit wafer with make-span is observed. As the rate of energy decreases with an increase in throughput, we observe this trend. Moreover, it can clearly be seen the energy per unit wafer in the current model is less than the baseline model. This is because even though the baseline model gives higher throughput, its energy consumption is far more than in the current model. Notably, the graph exhibits regions where the energy savings from machine on/off optimisation are especially pronounced, suggesting scenarios where machines might be running inefficiently or superfluously in a typical fab setting. However, it's crucial to acknowledge that the energy savings, while significant, may come with trade-offs in terms of machine wear, potential start-up delays, and other operational intricacies. This could lead to higher long-term capital and maintenance costs. Additionally, while energy savings are evident, there might be challenges in synchronisation, especially in a highly interconnected production line. Machines that are frequently switched off might not be immediately available when needed, leading to potential delays or bottlenecks. Hence, moving forward, it would be prudent to designate a maximum limit for machine restarts to prevent breakdowns.

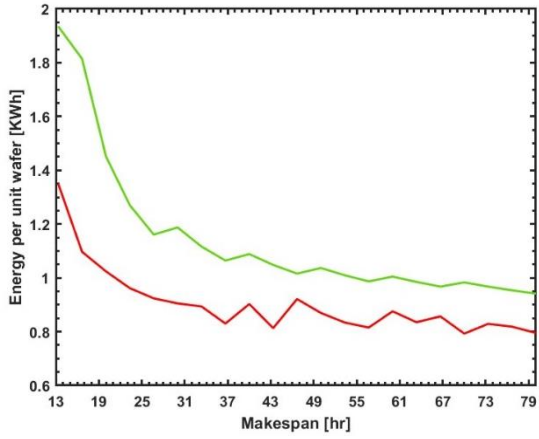

**Figure 11.** Comparing energy consumed per unit wafer for the baseline model (green solid line) and dynamic machine switching model (red solid line).

### 5.4 | First In First Out – FIFO dispatch strategy

This section highlights the importance of optimal scheduling on a shop floor to optimise energy consumption. The First in first out (FIFO) strategy ensures the wafer starting its operation first should be the one completing all the operations first as well. Here, we take a heuristic-based approach to find an optimal schedule, allowing the solver a flexible time horizon. We do not assign any objective function; an additional equation is introduced to the model discussed in section 3.3 to ensure the start and end times of the first and last operation, respectively, are in order.

$$\sum_{t'=t}^{t'=t} X_{k,j,m,t'} \leq \sum_{t'}^{t'=t-p_{k,m}} X_{k,j-1,m,t'} \quad (26)$$

In the context of a FIFO strategy, the above equation ensures that products are processed in the sequence of their arrival. Specifically, if job $j$ arrives after job $j-1$, it will not commence its processing until job $j-1$ has been allocated adequate time to complete its operation, even if multiple machines can perform the operation type $k$. This constraint upholds the core principle of FIFO: products introduced earlier should be processed ahead of those arriving later. This not only maintains a consistent processing order in line with their arrival sequence but also enhances predictability in the production environment. Predictable processing allows for better inventory management, timely response to disruptions, and streamlined production flows. Additionally, the energy analysis model introduced in section 5 is used with a constant throughput of 37 and allocating a flexible time horizon. The maximum time horizon is varied from 80 hours to 105 hours. The multi-objective formulation is ignored, as the aim of this test is to observe machine choices and energy consumption based on the FIFO heuristic.

The trend of energy consumption with make-span is displayed in Figure 12. It must be noted that the throughput is constant at 37 for each make-span value. As we do not assign any objective/constraint on make-span, the solver utilises the entire time horizon to complete the operations. The energy consumption corresponding to a minimum make-span of 81 hours is 23.10 KWh, which is consistent with the results discussed in the previous section. The instance with the maximum energy consumed is 31.82 KWh. As the make-span increases, we observe fluctuations in energy consumption. This suggests that the longer the tasks and processes are spread out, the more variability there is in how machines and resources are utilised, leading to different energy consumption patterns.

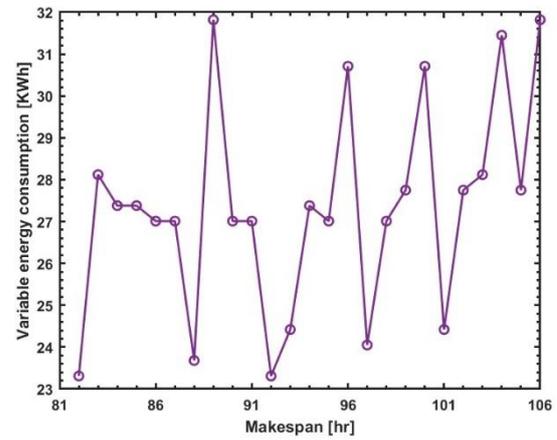

**Figure 12.** Energy consumption with different make-span and fixed throughput implementing FIFO heuristics.

The energy consumption might seem to follow some sort of periodic trend; essentially, it is the manifestation of choosing different (consistent) machine combinations by the solver. Certain make-span values seem to correspond to lower energy consumption, suggesting more efficient use of resources during those periods. Conversely, higher energy consumption at other make-span values indicates periods of inefficiency, perhaps due to machines being left idle or sub-optimal sequencing of tasks. This plot underscores the importance of scheduling not just for performance (make-span) but also for energy efficiency, as discussed in previous sections. The sequencing in FIFO, while straightforward, doesn't always guarantee the most efficient use of resources, as it doesn't prioritise tasks based on machine energy profiles or other optimisation parameters.

# 6 | CONCLUSION

This study elucidates two important questions prevalent in semiconductor manufacturing- first, how does energy consumption trade-off with throughput over a given time horizon? and second, how essential is optimal scheduling in order to reduce energy expenditure? To answer the first question, we formulate a flexible job shop scheduling problem using the six step Minifab dataset and state task network. The model is modified by assigning power rating (weights) to each machine and a multi-objective formulation is made to minimise energy consumption and maximise the throughput. The $\varepsilon-$ constraint method allows us to construct trade-off curves between energy and throughput. A linear increase of energy with throughput is observed. Although the gradient of energy decreases for a throughput of 30-37. This points towards the fact that the energy consumption per unit wafer decreases with an increase in make-span. A test case for fixed throughput is run, and the results give us distinct scenarios of energy consumption within bounds based on distinct machines selected. This gives us the flexibility to choose among available machines as some machines might be unavailable at certain times. Hence, based on machine selection, we can account for the energy expenditure.

To mimic a more realistic shop floor scenario, machines are dynamically switched on/off and account for a 'startup time' every time it is switched on. This results in a decrease of 2 wafers in throughput, but the energy consumption decreases by 7 KWh compared to the baseline model. Hence, machines with higher power ratings are utilised restrictively. Even the energy usage per unit wafer is significantly lower, hence pointing towards the major benefit of switching off machines when not required. The energy savings realized, especially in the context of large-scale manufacturing setups, could translate to significant financial and environmental benefits.

The second question can be answered by the FIFO heuristic model, where we do not assign an objective to minimise the make-span. The solver chooses to utilise the entire time horizon to fabricate 37 wafers, and the energy consumption is distinct based on machine choices and machine idle times. The maximum energy consumed is 31.82 KWh, which is a 36.46% increase than the optimised energy consumption.

As part of future research, it would be interesting to study the effect of introducing certain 'priority lots' (hot lots) in the scheduling model. Machine utilisation patterns and energy-saving strategies will be affected as the solver will prioritize these wafers before normal lots. This can be implemented in two ways: first, by assigning certain weights to wafers and amending the objective function to incentivise processing these wafers first. The second method can be introducing additional constraints for starting time for the first operation and ending time for the last operation for hot lots to be less than normal lots. Although this method will increase the size of the model significantly and would be computationally more expensive. Additionally, the current model could be refined further by increasing the size of lots and the number of machines (tools). The concept of buffers and batching should be implemented, and the model should be simulated for real fab datasets rather than test values to gain a deeper understanding of machine utilization and energy efficiency. These changes will be implemented and assessed in a companion paper in the future.

## ACKNOWLEDGEMENTS

The authors would like to extend their sincere gratitude to the Department of Chemical Engineering at Imperial College London for granting access to their computational resources and required licenses that were used in this research. We also acknowledge the collaborative partnership between Flexciton and Imperial College London, which significantly contributed to the progression of this work.

## STATEMENT FOR CONFLICT OF INTEREST

The authors declare no conflict of interest.